\input amstex
\documentstyle{amsppt}
\magnification=\magstep1
\hsize=5.2in
\vsize=6.8in
\centerline {\bf ON THE FUNDAMENTAL GROUP OF TYPE II$_1$ FACTORS} 
\vskip .1in
\centerline {\rm by}
\vskip .1in
\centerline {\rm SORIN POPA\footnote"*"{Supported by 
NSF-Grant 0100883}}

\address Math Dept  
UCLA, Los Angeles, CA 90095-1555\endaddress
\email popa\@math.ucla.edu\endemail

\vskip .2in
\centerline{\bf Dedicated to Richard V. Kadison}
\vskip .2in

\topmatter
\abstract We present here a shorter version of the 
proof of a result in ([Po2]), 
showing that  
the von Neumann factor associated with the  
group $\Bbb Z^2 \rtimes SL(2, \Bbb Z)$ has trivial fundamental group. 

\endabstract
\endtopmatter

The fundamental group of a type II$_1$ factor $M$ is the set 
of numbers $t>0$ for which the ``amplification'' 
of $M$ by $t$ is isomorphic to $M$, 
$\Cal F(M)=\{t>0 \mid M^t \simeq M\}$. 

This invariant for II$_1$ factors 
was considered by Murray and von Neumann 
in their pioneering paper ``Rings of operators IV'' ([MvN]),  
in connection with their notion of 
continuous dimension. They noticed that $\Cal F(M)=\Bbb R_+^*$ 
when $M$ is isomorphic to the approximately finite dimensional 
(or hypefinite) type II$_1$ factor, and 
more generally when it splits off such a factor. In  ([C1]) 
Connes showed that pairs $(M, t)$ with $t\in \Cal F(M), t\neq 1$, 
are at the ``core'' of the structure 
and classification of factors of type III. He also showed in ([C2]) 
that $\Cal F(M)$ reflects the rigidity properties 
of $M$, being countable whenever $M\simeq L(G_0)$ 
for $G_0$ a group with the property T of Kazhdan. A key 
application in Voiculescu's free probability theory 
was to show that $\Cal F(M)=\Bbb R_+^*$ for $M\simeq L(\Bbb F_\infty)$ 
([V], [Ra]), while the problem of calculating $\Cal F(M)$ 
for $M \simeq L(\Bbb F_n), 2\leq n < \infty$,  
remained its central open problem ([VDN], [Ra], [Dy]). 

On  the other hand, in the ergodic theory of groups, Gaboriau 
defined and calculated several invariants 
for orbit equivalence relations $\Cal R=\Cal R_{\sigma, \Gamma_0}$ coming 
from measure preserving 
actions $\sigma$ of groups $\Gamma_0$ on the probability space $(X, \mu)$,   
such as cost, $\ell^2$-Betti numbers, etc ([G1,2]). This  
allowed him to deduce that for a large class of groups $\Gamma_0$, including 
the free groups with finitely many generators $\Bbb F_n$, 
the fundamental group of $\Cal R$, 
$\Cal F(\Cal R)=\{t>0\mid \Cal R^t \simeq \Cal R\}$, is trivial. 
Equivalently, if $A=L^\infty(X, \mu) \subset M=A \rtimes \Bbb F_n$, 
is the II$_1$ factor with its 
Cartan subalgebra associated to the action of $\Bbb F_n$ ([FM]), then   
the fundamental group of the inclusion $A\subset M$ ([Po1]), 
$\Cal F(A\subset M)=\{t>0\mid \exists \theta : M^t \simeq M, \theta(A^t)=A\}$, 
is trivial.  

In ([Po2]), we introduced and studied the class 
of type II$_1$ factors $M$ 
having Cartan subalgebras $A \subset M$ with respect 
to which $M$ satisfies a combination of rigidity and weak-amenability 
properties, that we called HT. The prototype examples 
of HT factors in ([Po2]) 
are the cross-product algebras $M=L^\infty(\Bbb T^2, \mu) 
\rtimes_{\sigma_0} \Gamma_0$, with $\Gamma_0 \subset SL(2, \Bbb Z^2)$ 
subgroups of finite index 
(e.g., $\Gamma_0 = \Bbb F_n$) and $\sigma_0$ the restriction to $\Gamma_0$ 
of the action of $SL(2, \Bbb Z)$ on $\Bbb T^2 = \hat{\Bbb Z^2}$.  
The key result in ([Po2]) was a ``unique cross-product 
decomposition'' for HT factors, which implied that all 
invariants for $A \subset M$ (notably the 
$\ell^2$-Betti numbers in [G2]) are invariants 
for the HT factors $M$. As a consequence, we obtained  
many examples of type II$_1$ factors $M$ 
with trivial fundamental group, $\Cal F(M)=\{1\}$, 
including $M=L(\Bbb Z^2 \rtimes SL(2, \Bbb Z))=L^\infty(\Bbb T^2, \mu) 
\rtimes SL(2, \Bbb Z)$. In particular, this 
solved a long standing problem of Kadison, by showing that there 
exist factors $M$ with $M_{n\times n}(M) \not\simeq M$ 
for any $n\geq 2$ ([K1]; see also [Sa]). 

In this paper we give a short proof of the fact that $\Cal F(M)=\{1\}$, 
for many of the concrete examples of HT factors in ([Po2]).   
The proof is in fact not new, it is just the adaption 
to these particular cases of the proof in ([Po2]). Thus, by 
considering only algebras and subalgebras associated to a 
group-subgroup situation, we bypass the lengthy technicalities 
involved in developing the abstract property HT. 

The type II$_1$ 
factors $M$ that we consider here are defined as follows: 

Let $\sigma_0$ be 
the free, ergodic action of $SL(2, \Bbb Z)$ on $\Bbb T^2$, 
implemented by its action on $\Bbb Z^2$ 
(thus $\sigma_0$ preserves the Lebesgue measure 
$\mu_0$ on $\Bbb T^2$). Let $\Gamma_0 \subset SL(2, \Bbb Z)$ be 
any subgroup with finite index. Let $\sigma_1$ be any 
ergodic, measure preserving action of $\Gamma_0$ 
on some probability space $(X_1, \mu_1)$. Let 
$(X, \mu) = (\Bbb T^2, \mu_0)\times (X_1, \mu_1)$ and $\sigma = 
\sigma_0\times \sigma_1$ the corresponding product action 
of $\Gamma_0$ on $(X, \mu)$. We define 
$M$ to be the 
cross-product (or group-measure space) 
von Neumann algebras $M=L^\infty(X, \mu)
\rtimes_\sigma \Gamma_0$. The action $\sigma_0$ is  
well 
known to be mixing on $(\Bbb T^2, \mu_0)$. Since $\sigma_1$ 
is ergodic, $\sigma$ follows ergodic and thus  
$M$ follow type II$_1$ factors (see e.g., [D2], [Sa] or [T]). 

\proclaim{Theorem} The type II$_1$ factors $M$ constructed above  
have trivial fundamental group, $\Cal F(M)=\{1\}$. 
\endproclaim 

\proclaim{Corollary} For each finite $n \geq 2$ 
there exist free, ergodic, 
measure-preserving actions $\sigma$ of $\Bbb F_n$ on 
the probability space $(X, \mu)$ such that 
$M = L^\infty(X, \mu) \rtimes_\sigma \Bbb F_n$ has 
trivial fundamental group, $\Cal F(M)=\{1\}$. Moreover, the action $\sigma$ can 
be taken strongly ergodic, thus getting $M$ to be 
non-$\Gamma$, or non-strongly ergodic, 
thus getting $M$ to have the property $\Gamma$.
\endproclaim 

There are two crucial properties that enable us to prove 
the Theorem: 
on the one hand Kazhdan's 
rigidity of the inclusion 
$\Bbb Z^2 \subset \Bbb Z^2 \rtimes SL(2, \Bbb Z)$ ([Kaz]), 
which makes the entire subalgebra $A=L^\infty(X, \mu)$ 
sit ``somewhat rigidly'' inside $M=L^\infty(X, \mu)
\rtimes_\sigma \Gamma_0$; on the other hand, 
Haagerup's compact approximation property of the groups $\Gamma_0$ ([H]), 
which make the factors $pMp$ be in some sense 
``weakly amenable'' relative to  $Ap$, 
for any non-zero projection $p\in \Cal P(A)$. 

As a consequence of these properties, if there would exist 
an isomorphism $\theta$ of $pMp$ onto $M$, 
for some $p\in \Cal P(A), p\neq 0,1,$ then $M$ 
would have maximal abelian subalgebras $A$ and $\theta(Ap)$ 
with $A \subset M$ rigid and $M$ weakly 
amenable relative to $\theta(Ap)$. The 
heart of the argument consists in 
proving that if a type II$_1$ factor $M$ has 
two such maximal abelian subalgebras, $A$ and $\theta(Ap)$, 
then these subalgebras are conjugate by a unitary element in $M$ 
(this is what we earlier called ``unique 
decomposition result'' for $M$). But if this is the case, then 
the inclusions $A\subset M$ and $Ap \subset pMp$ follow isomorphic 
(via Ad$u\circ \theta$), 
thus implementing isomorphic measurable equivalence 
relations on the standard probability space ([FM]). When 
combined with the 
results in ([G1]) or ([G2]), this gives a contradiction. 

The split the details of this argument into several lemmas: 

\proclaim{Lemma 1} Let $A=L^\infty(X, \mu)  \subset L^\infty(X, \mu)
\rtimes_\sigma \Gamma_0=M$ be as defined above 
and let $\{u_g\}_{g\in \Gamma_0}\subset M$ be the canonical 
unitaries implementing the action $\sigma$. Let $p\in A$ 
be a non-zero projection and $\theta: pMp \rightarrow M$ 
an isomorphism. Denote $B=\theta(Ap)$ and $v_g = \theta(pu_gp), 
g\in \Gamma_0$. Then we have: 

$(a)$. $v_g$ are partial isometries in the normalizing groupoid 
$\Cal G\Cal N_M(B)$ of $B$ in $M$, i.e., $v_gv_g^* \in \Cal P(B)$, 
$v_g^*v_g \in \Cal P(B)$ and $v_gBv_g^* = Bv_gv_g^*$, 
$\forall g\in \Gamma_0$. Moreover, the 
spaces $v_gB, g\in \Gamma_0,$ are mutually orthogonal in $L^2(M)$ and we have   
$\overline{\Sigma_g v_gB} = M$. 

$(b)$. If $\varphi_0$ is a positive 
definite function on $\Gamma_0$ then 
the application $\Sigma_g v_g b_g \mapsto  
\Sigma_g \varphi_0(g) v_g b_g$, defined on finite sums with $b_g\in B$, 
extends to a unique unital, trace-preserving, completely 
positive, $B$-bilinear map $\phi_{\varphi_0}: 
M\rightarrow M$. Similarly, $\Sigma_g \hat{v_g b_g} \mapsto  
\Sigma_g \varphi_0(g) 
\hat{v_g b_g}$ extends to a unique bounded linear 
operator $T_{\varphi_0}$ on $L^2(M)$.  Moreover, we have  
$T_{\varphi_0}=\Sigma_g \varphi_0(g) v_ge_Bv_g^* = 
\Sigma_g \varphi_0(g) proj_{L^2(v_gB)}$, where 
$e_B= proj_{L^2(B)}$. Also, $T_{\varphi_0}(\hat{x}) 
=\hat{\phi_{\varphi_0}(x)}$, $\forall x\in M$. 
\endproclaim
\noindent 
{\it Proof}. $(a)$. Since $pu_gp= u_g \sigma_g^{-1}(p)p 
\in \Cal G\Cal N_{pMp}(Ap)$, we have $v_g=\theta(pu_gp) \in 
\Cal G\Cal N_M(B)$. Also, since $\Sigma_g u_gA$ 
is dense in $M$, $\Sigma_g pu_gp Ap= p(\Sigma_g u_g A)p$ 
is dense in $pMp$, so $\Sigma_g v_g B = \theta(\Sigma_g pu_gp Ap)$ 
is dense in $M$. Since $p(u_gA)p=u_g \sigma_{g^-1}(p)pA$ are 
mutually orthogonal, $v_gB, g\in \Gamma_0$, follow mutually 
orthogonal as well. 

$(b)$. The same argument as (1.1 in [H]) shows that 
$\phi_0(\Sigma u_g a_g) = \Sigma_g \varphi_0(g) u_g a_g$, 
with $a_g \in A$ all but finitely many equal to $0$, extends to a  
unique unital, trace-preserving, completely 
positive, $A$-bilinear map $\phi_0$ on $M$. Thus, 
$\phi'_0=\phi_0(p \cdot p)$ is a unital, trace-preserving, completely 
positive, $Ap$-bilinear map on $pMp$ and 
$\phi_{\varphi_0} = \theta \circ \phi_0' \circ \theta^{-1}$ 
is a unital, trace-preserving, completely 
positive, $B$-bilinear map on $M$. 

Since the vector spaces $v_g B, g\in \Gamma_0,$ are 
mutually orthogonal in $L^2(M)$, it follows 
that $\Sigma_g \varphi_0(g) v_ge_Bv_g^* = 
\Sigma_g \varphi_0 proj_{L^2(v_gB)}$ is a well defined 
positive operator on $L^2(M)$ and by the definitions it is trivial 
to see that it coincides with $T_{\varphi_0}$. 
\hfill Q.E.D. 

For the next lemma, 
recall that by a well known 
result of Kazhdan ([Kaz]) the inclusion of groups 
$G=\Bbb Z^2 \subset G_0=\Bbb Z^2 \rtimes SL(2, \Bbb Z)$ satisfies the 
following rigidity property: ``If a unitary representation of $G_0$ 
weakly contains the trivial representation of $G$ then it 
contains the trivial representation of $G$''. It is immediate 
to see that if we take $\Gamma_0 \subset SL(2, \Bbb Z)$ 
a subgroup of finite index, then $G=\Bbb Z^2 
\subset G_0=\Bbb Z^2 \rtimes \Gamma_0$ still satisfies 
the above rigidity condition. Also, 
the same proof as in ([AW]) shows that 
for a inclusions of groups $G \subset G_0$ with $G$ normal in $G_0$ 
this rigidity condition 
is equivalent to: ``If $\varphi_n$ 
are positive definite functions on $G_0$ with 
$\underset n  \to \lim\varphi_n(g)=1, \forall g\in \Gamma_0,$ 
then $\varphi_n$ converge to $1$ uniformly on $G$''. This 
in turn is clearly equivalent to the following: ``$\forall 
\varepsilon > 0$, $\exists F_0(\varepsilon)\subset G_0$ finite, 
$\delta_0(\varepsilon) >0$ such that if $\varphi$ is a positive 
definite function on $G_0$ with  
$|\varphi(g_0)-1|\leq \delta_0(\varepsilon), 
\forall g_0\in F_0(\varepsilon)$ then 
$|\varphi(h)-1|\leq \varepsilon, \forall h\in G$.'' (see [dHVa], [Jo] for 
more on this ``relative property T'' condition, 
which for arbitrary inclusions $G\subset G_0$ was first emphasized by 
Margulis.)

The Lemma below uses group rigidity in the same spirit as ([C], [CJ]): 
 
\proclaim{Lemma 2} Let $M$ and 
$G_0 = \Bbb Z^2 \rtimes \Gamma_0$ be as before. 
For any $\varepsilon > 0$ there exist 
a finite set $F(\varepsilon) \subset G_0$ and 
$\delta(\varepsilon)>0$ such that if 
$\phi : M \rightarrow M$ is a unital 
completely positive map with $\|\phi(u_{g_0})-u_{g_0}\|_2 \leq 
\delta(\varepsilon)$, $\forall g_0 \in F(\varepsilon)$, 
then $\|\phi(u_h)-u_h\|_2 \leq 
\varepsilon$, $\forall h\in \Bbb Z^2$. 
\endproclaim
\noindent
{\it Proof}. Let $\phi$ be a unital completely positive map 
satisfying $\|\phi(u_{g_0})-u_{g_0}\|_2 \leq 
\delta_0(\varepsilon^2/2)$, $\forall g_0 \in F_0(\varepsilon^2/2)$. 
Define $\varphi(g)=\tau(\phi(u_g)u_g^*), g\in G_0$. Then 
$\varphi$ is positive definite and by the 
Cauchy-Schwartz inequality we have: 
$$
|\varphi(g_0)-1|=|\tau(\phi(u_{g_0})u_{g_0}^*) - \tau(u_{g_0}u_{g_0}^*)|
$$
$$ 
\leq \|\phi(u_{g_0})-u_{g_0}\|_2 \|u_{g_0}^*\|_2 \leq \delta_0(\varepsilon^2/2),
$$
for any $g_0 \in F_0(\varepsilon^2/2)$.   
Since by 
Kadison's inequality we have 
$\phi(u_h)\phi(u_h^*) \leq \phi(u_hu_h^*) = 1$, 
it follows that for all $g \in G_0$ we have 
$\|\phi(u_g)-u_g\|_2^2 = 1 - 2 {\text{\rm Re}} 
\tau(\phi(u_g)u_g^*) + \tau(\phi(u_g)\phi(u_g^*)) 
\leq 2(1 - {\text{\rm Re}} \varphi(g))$. Thus, if $h \in G$ we get: 
$$
\varepsilon^2/2 \geq |\varphi(h)-1| \geq (1 - 
{\text{\rm Re}} \varphi(h)) 
\geq \|\phi(u_h)-u_h\|_2^2/2. 
$$

Thus, if we define $F(\varepsilon)=F_0(\varepsilon^2/2)$, 
$\delta(\varepsilon)=\delta_0(\varepsilon^2/2)$, 
then we are done. 
\hfill Q.E.D. 

\proclaim{Lemma 3} Let $A, B \subset M$ be as before. There 
exists a finite subset $F\subset G_0$ such that 
the projection $f = \Sigma_{g\in F} v_ge_Bv_g^* = 
\Sigma_{g\in F} proj_{L^2(v_gB)}$ 
satisfies $\|f(\hat{u_h}) - \hat{u_h}\|_2 \leq 1/4$, $\forall h\in \Bbb Z^2$. 
\endproclaim 
\noindent 
{\it Proof}. Since $SL(2, \Bbb Z)$ 
contains free groups as subgroups with finite index, 
by ([H]) it has the Haagerup approximation 
property and 
so do all its subgroups 
of finite index $\Gamma_0$. Let $\varphi_n: \Gamma_0 \rightarrow \Bbb C$ 
be positive definite functions with $\varphi_n \in c_0(\Gamma_0)$, 
$\varphi_n(e)=1, \forall n,$  
and $\underset n  \to \lim\varphi_n(g)=1, \forall g\in \Gamma_0.$

By the definition of $\phi_{\varphi_n}=\phi_n$ 
in Lemma 1, it follows that $\phi_n$ tend 
to $id_M$ in the 
point-norm $\|\quad\|_2$-topology.  Let $n$ be large enough 
so that $\|\phi_n(u_{g_0})-u_{g_0}\|_2 \leq \delta(1/16)$, 
$\forall g_0\in F(1/16)$, where $F(1/16)\subset G_0, 
\delta(1/16)$ are 
as given by Lemma 2. It follows that $\|\phi_n(u_h)-u_h\|_2 \leq 
1/16$, $\forall h\in \Bbb Z^2$. Thus, if we let $T=T_{\varphi_n}$ 
then $0\leq T \leq 1$ and $\|T(\hat{u_h}) - \hat{u_h}\|_2 \leq 1/16$, 
$\forall h\in \Bbb Z^2$. Let $f$ be the spectral 
projection of $T$ corresponding to [$1/16, 1$]. Thus, 
if $F$ is the set of all $g \in \Gamma_0$ with 
$\varphi_n(g) \geq 1/16$ then $F$ is finite and 
$f = \Sigma_{g\in F} v_ge_Bv_g^* = 
\Sigma_{g\in F} proj_{L^2(v_gB)}$. Moreover, we have 
$\|f(\hat{u_h}) - u_h\|_2 \leq 1/4$, $\forall h\in \Bbb Z^2$. 
Indeed, for if we would have $\|(1-f)(\hat{u_h})\|_2 > 1/4$ 
for some $h\in \Bbb Z^2$, then 
$$
1/16 \geq \|(1-f)(T(\hat{u_h}) - \hat{u_h})\|_2 
$$
$$
\geq \|(1-f)(\hat{u_h})\|_2 - \|(1-f)T(\hat{u_h})\|_2 > 1/4 - 1/16=1/16, 
$$
a contradiction. 
\hfill Q.E.D. 

\proclaim{Lemma 4} Let $A_0 \subset M$ be the 
von Neumann algebra generated by $\{u_h\}_{h \in \Bbb Z^2}\subset M$. 
There exists $\xi \in L^2(M), \xi \neq 0,$ such that 
$A_0 L^2(\xi B) \subset L^2(\xi B)$. 
\endproclaim
\noindent 
{\it Proof}. To obtain from 
the projection $f$ in Lemma 3 a left $A_0$ module 
of the form $L^2(\xi B)$, we first use a trick 
inspired from ([Chr]). Thus, 
let $\langle M, B\rangle$ be the von Neumann 
algebra generated in $\Cal B(L^2(M))$ 
by $M$ (regarded as the algebra 
left multiplication operators by elements in $M$) and by the 
orthogonal projection $e_B$ of $L^2(M)$ onto $L^2(B)$. Note  
that if $x\in M$ then $e_Bxe_B = E_B(x)e_B$ 
and that $\langle M, B \rangle$ is the weak closure of the 
$*$-algebra sp$\{xe_B y \mid x,y \in M\}$. Also, 
$\langle M, B \rangle$ coincides with the commutant 
in $\Cal B(L^2(M))$ of the operators 
of right multiplication by elements in $B$, i.e., with 
$J_M B' J_M$, 
where $J_M$ is the canonical conjugation on $L^2(M)$, 
defined by $J_M(\hat{x})=\hat{x^*}, x\in M$. (See [Ch] or [J] 
for details on $\langle M, B \rangle$ and on all this 
construction, called the {\it basic construction}). 

Since $e_B$ has 
central support 1 in $\langle M, B\rangle$ and $e_B 
\langle M, B \rangle e_B = B e_B$, it follows that 
$\langle M, B \rangle$ is of type I. Moreover, $Tr(xe_B y) = \tau(xy), 
x, y \in M$ defines a unique normal faithful semifinite trace on 
$\langle M, B \rangle$. 

Let $f = \Sigma_{g\in F} v_ge_Bv_g^*$ be the projection 
given by Lemma 3. Note that $f \in \langle M, B \rangle$ and 
$Tr(f) = \Sigma_{g\in F} \tau(v_gv_g^*) < \infty$.  
We denote by 
$K = {\overline{\text{\rm co}}}^{\text{\rm w}} 
\{u_h f u_h^* \mid h \in \Bbb Z^2 \}$. Since $K$ is weakly closed and  
$0\leq a \leq 1, \forall a\in K$, $K$ is weakly 
compact. Moreover, $Tr(a) \leq Tr(f), \forall a\in K$,  
so $K$ is also contained in the Hilbert space 
$L^2(\langle N, B \rangle, Tr)$, where it is still  
weakly compact. 

Let $k\in K$ be the unique element of minimal 
norm $\| \quad \|_{2,Tr}$ in $K$. 
Since $u_h K u_h^* = K$ and 
$\|u_hku_h^*\|_{2, Tr} = \|k\|_{2, Tr}, \forall h\in \Bbb Z^2$, 
by the uniqueness of $k$ 
it follows that $u_hku_h^*=k, \forall h\in \Bbb Z^2$.  
Thus, $k\in A_0'\cap 
\langle M, B \rangle$, $0 \leq k \leq 1$ and $Tr(k) \leq 1$.  

Since 
$$
Tr(e_Bu^*_hfu_h) = 
Tr(fu_he_Bu_h^*f)
$$
$$
= Tr ((\Sigma_{g\in F} u_ge_Bu_g^*) 
u_h e_Bu^*_h(\Sigma_{g'\in F} u_{g'}e_B u_{g'}^*)) 
$$
$$
=\Sigma_{g,g'\in F} \tau(u_gE_B(u_g^*u_h)E_B(u_h^*u_{g'})u_{g'}^*) 
$$
$$
= \|f(\hat{u_h})\|_2^2 \geq (3/4)^2, 
$$
it follows that $Tr(e_Ba) \geq 9/16, \forall a\in K$. 
In particular, $Tr(e_Bk) \geq 9/16$, so $k\neq 0$. 
Thus, if $e$ is a spectral projection 
of $k$ corresponding to $(c, \infty)$ 
for some appropriate $c > 0$, then 
$e \in A_0'\cap \langle M, B \rangle$, $e\neq 0$ and  
$Tr(e) < \infty$. Thus, $e\langle M, B \rangle e$ 
is a finite type I von Neumann algebra. Let $A_1 \subset 
e\langle M, B \rangle e$ be a maximal abelian subalgebra 
containing $A_0e$. By ([K2]) $A_1$ contains an abelian 
projection $e_0\neq 0$ of $e\langle M, B \rangle e$. Thus, 
$e_0$ is an abelian projection in $\langle M, B \rangle$. 
By von Neumann's theorem, there exists $\xi L^2(M)$ 
such that $e_0$ is the orthogonal projection onto $
\overline{\xi B}=L^2(\xi B)$. Since $e_0\in A_1$ 
commutes with $A_0$, we have $A_0 L^2(\xi B) \subset L^2(\xi B)$. 
\hfill Q.E.D. 

\proclaim{Lemma 5}  There exists a non-zero partial isometry $w\in 
M$ such that $w^*w \in A, ww^* \in B$ and 
$wAw^* \subset Bww^*$. 
\endproclaim
\noindent 
{\it Proof}. Let $A_0$ be  
as in Lemma 4. By construction, we readily see that $A_0'\cap M=A$. 
Let $\xi \in L^2(M), \xi \neq 0,$ be so 
that $A_0 \xi \subset L^2(\xi B)$ (cf. Lemma 4). 
By regarding $\xi$ as a 
square summable operator affiliated 
with $M$ and by replacing it by $\xi E_B((\xi^*\xi)^{-1/2})$, 
we may assume $q_0=E_B(\xi^*\xi)$ is a projection in $B$. 
Note that we then have $L^2(\xi B) = \xi L^2(B)$. Let 
$p_0\in A_0$ be the minimal projection in $A_0$ 
with the property that $(1-p_0)\xi=0$. We denote 
$\psi(a_0)=E_B(\xi^*a_0\xi), a_0\in A_0p_0$, and note that $\psi$ is a 
unital, normal, faithful, 
completely positive map from $A_0p_0$ into $Bq_0$. 

Also, since $a_0 \xi \in L^2(\xi B) = \xi L^2(Bq_0)$, it follows that 
$a_0 \xi = \xi \psi(a_0), \forall a_0 \in A_0$. Indeed, 
for if $a_0 \xi = \xi b$, for some $b \in L^2(B)q_0$, then 
$\xi^*a_0\xi = \xi^*\xi b$ and 
so 
$$
\psi(a_0)= E_B(\xi^*a_0\xi) = E_B(\xi^*\xi b)= E_B(\xi^*\xi)b=b. 
$$
Thus, 
for $a_1, a_2 \in A_0$ we get $a_1 a_2 \xi = a_1 \xi \psi(a_2) 
= \xi \psi(a_1) \psi(a_2)$. Since we also have 
$(a_1a_2) \xi = \xi \psi(a_1a_2)$, this shows that 
$\psi(a_1a_2) = \psi(a_1)\psi(a_2)$. Thus, $\psi$ is a 
$*$-morphism of $A_0$ into $A$. 

Thus, since $a_0 \xi = \xi \psi(a_0), \forall a_0 \in A_0,$ 
and $\psi$ is a $*$-isomorphism, 
by a standard trick, if $v\in M$ is the 
partial isometry in the polar decomposition of $\xi$ 
with the property that the right supports of $\xi$ and $v$ coincide,
then $p=vv^* \in A_0p_0'\cap p_0Mp_0=Ap_0$, 
$q=v^*v \in \psi(A_0)'\cap q_0Mq_0$ and 
$a_0 v = v \psi(a_0), \forall a_0 \in A_0$. Since 
$$
q(\psi(A_0)'\cap q_0Mq_0)q=(\psi(A_0)q)'\cap qMq
$$
$$
= (v^* A_0 v)'\cap 
qMq = v^*(A_0q_0'\cap q_0Mq_0)v = v^*Av, 
$$ 
it follows that $q$ is an abelian projection in $\psi(A_0)'\cap q_0Mq_0$. 
But $Bq_0$ is maximal abelian in $q_0Mq_0$, 
so by ([K2]) it follows that there 
exists a projection $q'\in Bq_0$ 
such that $q$ and $q'$ are equivalent in $\psi(A_0)'\cap q_0Mq_0$, 
say via a partial isometry $v' 
\in \psi(A_0)'\cap q_0Mq_0$. Thus, if we put $w=(vv')^*$ then 
$w^*w \in A$, $ww^* \in B$ and $wAw^* \subset B$. 
\hfill Q.E.D.

\proclaim{Lemma 6} There exists a unitary element $u\in M$ 
such that $uAu^* = B$. 
\endproclaim 
\noindent 
{\it Proof}. Note first that both $A$ and $B$ are 
regular maximal abelian subalgebras in $M$ in the sense 
of ([D1]; N.B.: Such algebras are called Cartan subalgebras in [FM]). 
Indeed, this is because $\Sigma_g u_g A$ 
(resp. $\Sigma_g v_g B$) is dense in $M$, implying 
that $\Cal G\Cal N(A)$ (resp. $\Cal G\Cal N(B)$) generates the 
von Neumann algebra $M$, so ([Dy]) applies.  

Let $w \in M$ be the non-zero partial isometry provided 
by Lemma 5. By cutting  
$w$ from the right with a smaller projection in $A$, we may 
clearly assume $\tau(ww^*) = 1/n$ for some integer $n$. 

Since $A, B$ are regular in $M$, by ([Dy]) there exist 
partial isometries $v_1, v_2, ..., v_n$, 
respectively $w_1, w_2, ..., w_n$ 
in the normalizing groupoids of $A_0$ respectively $A_1$ such that 
$\Sigma_i v_iv_i^* = 1, \Sigma_j w_jw_j^* = 1$ and 
$v_i^*v_i = w^*w, w_j^*w_j = ww^*, \forall i,j$. But then 
$u=\Sigma_i w_i w v_i^*$ is a unitary element and $vA_0v^*=B$.
\hfill Q.E.D. 

We summarize Lemmas 1-6 as the following ``unique 
decomposition result'':  

\proclaim{Proposition 7} If $p\in A$ 
is a non-zero projection and $\theta: pMp \rightarrow M$ is 
an isomorphism, then there exists a unitary element $u \in M$ 
such that ${\text{\rm Ad}}u(\theta(Ap))=A$. Thus, 
the inclusion $Ap\subset pMp$ and $A \subset M$ follow 
isomorphic. 
\endproclaim

\noindent
{\it Proof of the Theorem}. If $t \in \Cal F(M)$ and $t < 1$, then 
let $p\in M$ be a projection with trace $\tau(p) = t$ and 
$\theta : pMp \simeq M$ an isomorphism. Without 
loss of generality, we may take $p\in A$. By Proposition 7, 
$Ap \subset pMp$ and $A\subset M$ follow isomorphic. 
By ([FM]), this 
implies the countable, measure 
preserving equivalence relations $\Cal R_{Ap\subset pMp}$ and 
$\Cal R_{A\subset M}$ are isomorphic. But by ([G1] or [G2]), 
this is a contradiction. 
\hfill Q.E.D. 
\vskip .05in
\noindent
{\it Proof of the Corollary}. Since any free group with 
finitely many generators $\Bbb F_n, n\geq 2,$ can 
be embedded into $SL(2, \Bbb Z)$ with finite index, 
we can take $\Gamma_0=\Bbb F_n$ in the 
construction of $M$, proving the first part of the statement. 

If $M=L^\infty(X, \mu) \rtimes \Bbb F_n$, 
then any central sequence for $M$ is contained  
in $L^\infty(X, \mu)$ (see e.g., [Po3]). 
Thus, if we take $\sigma_1$ to be trivial, 
i.e., $X=X_0=\Bbb T^2$ with $\sigma=\sigma_0$, 
then by taking into account that $\sigma_0$ is strongly ergodic 
(cf. [Sc]), it follows that $M$ has no non-trivial central 
sequences, i.e., it is non-$\Gamma$. 

On the other hand, 
if we take $\sigma_1$ to be any free ergodic action 
of $\Bbb Z$ on a probability space $(X_1, \mu_1)$ (e.g., 
a Bernoulli shift action of $\Bbb Z$) composed with 
the quotient map $\Bbb F_n \rightarrow \Bbb Z$, 
then $\sigma=\sigma_0 \times \sigma_1$ gives an ergodic 
but not strongly ergodic action of $\Bbb F_n$ 
on $(X, \mu)$, thus getting $M$ to have 
property $\Gamma$. 
\hfill Q.E.D. 
\vskip .1in 
\noindent
{\bf Remark}. Note that the arguments in Lemmas 1-6 prove 
a more general ``unique decomposition result'' than the Proposition 7, 
namely: Let $G^i$ be commutative groups and $\Gamma^i_0$ 
some groups satisfying Haagerup's compact approximation property, 
$i=1,2$. Assume $\Gamma_0^i$ acts on $G^i$ mixingly and outerly, 
and so that the inclusion $G^i \subset G^i\rtimes \Gamma_0^i$ has 
the relative property T, for each $i=1,2$, in the sense of Margulis ([dHVa]; 
see [Va] for examples). For each $i=1,2$ let $\sigma_0^i$ 
be the corresponding trace preserving action of $\Gamma_0^i$ on $L(G_i)$ 
and let $\sigma_1^i$ be an ergodic measure-preserving 
action of $\Gamma_0^i$ on a probability space 
$(X_1^i, \mu_1^i)$ such that the action $\sigma^i = \sigma_0^i \otimes 
\sigma_1^i$ of $\Gamma_0^i$ on $A_i = L(G^i)\overline{\otimes} 
L^\infty(X_1^i, \mu_1^i)$ is ergodic. Let $M_i = A_i \rtimes_{\sigma^i} \Gamma_0^i$. If there exists some isomorphism $\theta$ from 
$p_1M_1p_1$ onto $p_2M_2p_2$ 
for some $p_i \in \Cal P(A_i)$ then there exists a 
unitary element $u\in p_2M_2p_2$ such that 
Ad$u \circ \theta (A_1p_1) = A_2p_2$. 
Note that this statement covers most of the concrete examples 
of HT factors in ([Po2]), but that it only proves the uniqueness 
of the ``concrete'' HT Cartan subalgebras. 

The proof of the general ``unique 
decomposition result'' for arbitrary HT factors in ([Po2]) requires 
a full discussion of the operator algebraic 
concept of (relative) 
property H and T. As a consequence, although going 
along the same lines, the proof there is unavoidably 
lenghtier. The full generality 
of that result is needed in order to show that the class of 
abstract HT factors 
is well behaved to elementary operations 
(amplification, tensor product, finite index extension/restriction), 
to make
the definition of Betti numbers for factors in the class HT be 
more conceptual (while relying on [G2]) 
and to show that they check the appropriate formulas. 
It is also needed in order  
to investigate the structure of finite Connes' 
correspondences and of subfactors of finite 
Jones index of  
HT factors (Sec. 7 in [Po2]), which in turn  
is of interest  because it relates a number theory framework to the theory of type II$_1$ factors. 

\head References\endhead

\item{[AW]} C. Akemann, M. Walters: 
{\it Unbounded negative definite functions}, 
Canad. J. Math. {\bf 33} (1981), no. 4, 862--871

\item{[Ch]}
E. Christensen: {\it Subalgebras of a finite algebra}, 
Math. Ann. {\bf 243} (1979), 17-29.

\item{[C1]} A. Connes: {\it Une classification 
des facteurs de type} III, Ann. \'Ec. Norm. Supe. 
{\bf 6} (1973), 133-252. 

\item{[C2]} A. Connes: {\it A type II$_1$ 
factor with countable fundamental group}, J. Operator 
Theory {\bf 4} (1980), 151-153. 

\item{[CJ]} A. Connes, V.F.R. Jones: {\it Property T 
for von Neumann algebras}, Bull. London Math. Soc. {\bf 17} (1985), 
57-62. 

\item{[D1]}
J. Dixmier: {\it Sous-anneaux ab\'{e}liens maximaux 
dans les facteurs de type fini}, Ann. Math. {\bf 59} (1954), 279-286. 

\item{[D2]} J. Dixmier: ``Les alg\`ebres d'op\'erateurs sur l'espace Hilbertien
(Alg\`ebres de von Neumann)'', Gauthier-Villars, Paris, 1957, 1969.

\item{[Dy]} H. Dye: {\it On groups of measure preserving transformations}, 
I, II, Amer. J. Math. {\bf 81} (1959), 119-159, and {\bf 85} (1963), 
551-576.

\item{[Dyk]} K. Dykema: {\it 
Interpolated free group factors}, Duke Math J. {\bf 69} (1993), 97-119.

\item{[FM]} J. Feldman, C.C. Moore: {\it Ergodic equivalence relations, 
cohomology, and von Neumann algebras I, II}, Trans. Amer. Math. 
Soc. {\bf 234} (1977), 289-324, 325-359. 

\item{[G1]} D. Gaboriau: {\it Cout des r\'elations d'\'equivalence 
et des groupes}, Invent. Math. {\bf 139} (2000), 41-98. 

\item{[G2]} D. Gaboriau: {\it Invariants $\ell^2$ de relations 
d'\'equivalence et de groupes}, preprint May 2001. 

\item{[H]} U. Haagerup: {\it An example of non-nuclear C$^*$-algebra 
which has the metric approximation property}, Invent. Math. 
{\bf 50} (1979), 279-293. 

\item{[dHVa]} P. de la Harpe, A. Valette: ``La propri\'et\'e T 
de Kazhdan pour les 
groupes localement compacts'', Ast\'erisque {\bf 175} (1989). 

\item{[Jo2]} P. Jolissaint: {\it On the relative property T}, 
preprint 2001.  

\item{[J]} V.F.R. Jones : {\it Index for subfactors}, Invent. Math. 
{\bf 72} (1983), 1-25.  

\item{[K1]} R.V. Kadison: {\it Problems on von Neumann algebras}, 
Baton Rouge Conference 1867, unpublished. 

\item{[K2]} R.V. Kadison: {\it Diagonalizing matrices}, Amer. Math. 
Journal (1984), 1451-1468. 

\item{[Kaz]} D. Kazhdan: {\it Connection of the dual space of a group 
with the structure of its closed subgroups}, Funct. Anal. and its Appl. 
{\bf1} (1967), 63-65. 

\item{[MvN]} F. Murray, J. von Neumann: {\it Rings of operators IV}, 
Ann. Math. {\bf44} (1943), 716-808.

\item{[Po1]} S. Popa: {\it Correspondences}, INCREST preprint 1986, 
unpublished. 

\item{[Po2]} S. Popa: {\it On a class of type} II$_1$ {\it factors with Betti numbers invariants}, MSRI preprint 2001/024, revised form 
August 2002. 

\item{[Po3]} S. Popa: {\it Maximal injective subalgebras in factors associated 
with free groups}, Adv. in Math., {\bf 50} (1983), 27-48. 

\item{[Ra]} F. Radulescu: {\it Random matrices, 
amalgamated free products and 
subfactors of the von {Neumann} algebra of a free group, 
of noninteger 
index}, Invent. Math. {\bf 115} (1994), 347--389.

\item{[Sa]} Sakai: ``C$^*$-algebras and W$^*$-algebras'', Springer-Verlag, 
Berlin-Heidelberg-New York, 1971.  

\item{[Sc]} K. Schmidt: {\it Asymptotically invariant 
sequences and an action of $SL(2, \Bbb Z)$ on the 
$2$-sphere}, Israel. J. Math. {\bf 37} (1980), 193-208.

\item{[T]} M. Takesaki: ``Theory 
of operator algebras I'', Springer, Berlin. 

\item{[Va]} A. Valette: {\it Semi-direct products with 
the property (T)}, preprint, 2001. 

\item{[V]} D. Voiculescu: {\it 
Circular and semicircular systems and free product 
factors}, ``Operator Algebras, 
Unitary Representations, Enveloping Algebras, 
and Invariant Theory'', 
Progress in Mathematics, {\bf 92}, Birkhauser, Boston, 
1990, pp. 45-60.

\item{[VDN]} 
D. Voiculescu, K. Dykema, A. Nica, ``Free random variables'', CRM 
monograph series, vol. 1, American Mathematical Society, 1992.

\bye